\documentclass[submission]{FPSAC2024}

\newtheorem{theorem}{Theorem}[section]
\newtheorem{thm-def}[theorem]{Theorem/Definition}
\newtheorem{lemma}[theorem]{Lemma}
\newtheorem{open}[theorem]{Open Problem}

\newtheorem{definition}[theorem]{Definition}
\newtheorem{corollary}[theorem]{Corollary}
\theoremstyle{definition}
\newtheorem{example}[theorem]{Example}
\theoremstyle{definition}

\theoremstyle{remark}

\DeclareMathOperator{\osp}{B}  

\DeclareMathOperator{\codim}{codim}

\newcommand{\given}{\, | \,}

\usepackage{setspace}

\newcommand{\Fullrk}{\mathcal{F}_{k \times n}(\mathbb{C})}
\newcommand{\W}{\mathcal{W}}

\usepackage{lipsum}

\title[Brewing Fubini-Bruhat Orders]{Brewing Fubini-Bruhat Orders}

\author[Billey and Ryan]{Sara C. Billey\addressmark{1}, \and Stark Ryan\addressmark{2}\thanks{\href{mailto:billey@math.washington.edu}{billey@math.washington.edu}. \href{mailto:drstarkryan@gmail.com}{drstarkryan@gmail.com}. The authors were partially supported by the Washington
  Research Foundation and  DMS-1764012.}}

\address{\addressmark{1}Department of Mathematics, University of Washington, Seattle, WA \\ \addressmark{2}Ghost Autonomy, Mountain View, CA}

\received{\today}

\abstract{ The Bruhat order on permutations arises out of the study of
Schubert varieties in Grassmannians and flag varieties, which have
been important for over 100 years \cite{BB,E,F,H,KL}. The purpose of
this paper is to study variations on this theme related to
subvarieties of the spanning line configurations $X_{n,k}$ as defined
by Pawlowski and Rhoades \cite{PR}.  These subvarieties are indexed by
Fubini words, or equivalently by ordered set partitions.  Three
natural partial orders arise in this context; we refer to them as the
decaf, medium roast, and espresso orders.  The decaf order is a
generalization of the weak order on permutations defined by covering
relations using simple transpositions.  The medium roast order is a
generalization of the (strong) Bruhat order defined by the closure
relationship on the subvarieties.  The espresso order is the
transitive closure of a relation based on intersecting subvarieties.
Many properties of Schubert varieties and Bruhat order extend to one
or more of the three Fubini-Bruhat orders.  We examine some of the
many possibilities in this work.  }


\usepackage[backend=bibtex]{biblatex}
\begin{document}

\maketitle

\section{Introduction}\label{sec:intro}

For positive integers $k \leq n$, a \textbf{Fubini word} $w=w_1 \cdots
w_n$ represents a surjective map $w:[n] \rightarrow [k]$. We denote a
Fubini word by its \textbf{one-line notation}, an ordered list
$w=w_1w_2 \cdots w_n$, where $w_i=w(i)$.  We denote by $\W_{n,k}$ the
Fubini words of length $n$ on the alphabet $[k]$.  For $k=n$, a Fubini
word $w \in \W_{n,n}$ is exactly a permutation in $S_n$, and the
one-line notation for $w$ is the same whether $w$ is viewed as a
Fubini word or a permutation.  The bijection between Fubini words and
ordered set partitions maps $w \in \W_{n,k}$ to
$\osp(w)=\osp_1|\osp_2|\ldots|\osp_k$ where $\osp_i =\{j \in [n] \, |
\, w_j =i \}$.  Hence the number of Fubini words in $\W_{n,k}$ is
$k!S(n,k)$ where $S(n,k)$ is the Stirling number of the second kind
\cite[A000670,A019538]{OEIS}.

Let $\Fullrk$ be the set of full rank $k \times n$ matrices with no
zero columns.  Such matrices have a Bruhat decomposition into orbits
\begin{equation}\label{eq:Full.Bruhat.Decomposition}
 \Fullrk = \bigsqcup_{w \in \W_{n,k}} B_{-}^{(k)} M_{w}
B_{+}(w) 
\end{equation}
where $M_{w}$ is the analog of a permutation matrix with a 1 in
position $(w_{j},j)$ and 0's elsewhere, $B_{-}$ and $B_{+}$ are the
set of invertible lower and upper triangular matrices respectively and
the superscript indicates their size, and $B_{+}(w) $ is the subgroup
of the $n \times n$ invertible upper triangular matrices $A$ such that
$M_{w}A \in \Fullrk$.  Every matrix in the double orbit $B_{-}^{(k)}
M_{w} B_{+}(w)$ can be written in many ways as a triple product, thus
it can be useful to chose canonical representatives.  Let
$U=U_{-}^{(k)}$ be the set of lower unitriangular matrices in
$GL_{k}(\mathbb{C})$, and let $T=T^{(n)}$ be the set of diagonal
matrices in $GL_{n}(\mathbb{C})$.  Pawlowski and Rhoades \cite{PR}
defined the \textbf{pattern matrices} $P_{w}$ indexed by words $w \in
\W_{n,k}$ to be a specific set of orbit representatives such that each
$M \in B_{-}^{(k)} M_{w} B_{+}(w) $ can be written uniquely as a
product $M=XYZ$ with $X \in U, Y \in P_{w}$, and $Z \in T$
\cite[Lem. 3.1 and Prop. 3.2]{PR}.  See \Cref{sec:background} for more
details.  Thus, we have an \textbf{efficient Bruhat decomposition}
\begin{equation}\label{eq:efficient.bruhat.decomp}
 \Fullrk = \bigsqcup_{w \in \W_{n,k}} U P_{w} T. 
\end{equation}


Under right multiplication, every $T$-orbit of $\Fullrk$ determines an
ordered list of $n$ 1-dimensional subspaces whose vector space sum is
$\mathbb{C}^k$ via its ordered list of columns.  The set of such
``lines'' in $\mathbb{C}^k$ is the $(k-1)$-dimensional complex
projective space $\mathbb{P}^{k-1}$.

\begin{definition}\label{def:spanning.line.config}\cite[Def. 1.3]{PR}
A \textbf{spanning line configuration} $l_{\bullet} =
(l_1,\ldots,l_n)$ is an ordered $n$-tuple in the product of projective
spaces $(\mathbb{P}^{k-1})^{n}$ whose vector space sum is
$\mathbb{C}^k$.  Let
\begin{equation}X_{n,k}= \Fullrk/T = \{l_{\bullet} = (l_1,\ldots,l_n) \in
(\mathbb{P}^{k-1})^n \hspace{0.05in} | \hspace{0.05in} l_1 + \cdots +
l_n = \mathbb{C}^k\} 
\end{equation}
be the \textbf{space of spanning line configurations} for $1 \leq
k\leq n$.
\end{definition}

In 2017, Pawlowski and Rhodes initiated the study of the space of
spanning line configurations \cite{PR}.  They observed and proved the
following remarkable properties.  The projection of $X_{n,n}=GL_{n}/T$
to the flag variety $GL_{n}/B_{+}^{(n)}$ is a homotopy equivalence, so
they have isomorphic cohomology rings.  More generally, $X_{n,k}$ is
an open subvariety of $(\mathbb{P}^{k-1})^n$, hence it is a smooth
complex manifold of dimension $n(k-1)$.  The cohomology ring of
$X_{n,k}$ may be presented as the ring
\[
R_{n,k} = \mathbf{Z}[x_{1},\dotsc , x_{n}]/\langle x_{1}^{k},\dotsc ,
x_{n}^{k}, e_{n-k+1},\dotsc, e_{n}\rangle
\]
defined by Haglund-Rhoades-Shimozono \cite{HRS}, generalizing the
coinvariant algebra and Borel's theorem $H^{*}(GL_{n}/B)\cong
R_{n,n}$.  Here, $e_{i}$ is the $i^{th}$ elementary symmetric function
in $x_{1},\dots , x_{n}$.  Furthermore, there is a natural $S_{n}$ action on
$n$-tuples of lines inducing an $S_{n}$ action on the cohomology ring
of $X_{n,k}$, which is isomorphic to $R_{n,k}$ as a graded
$S_{n}$-module.  See also \cite{Griffin} for another geometric
interpretation of $R_{n,k}$. The efficient Bruhat decomposition gives
rise to a cellular decomposition
\[
X_{n,k} = \bigsqcup_{w \in \W_{n,k}} U P_{w}.  
\]
Let $C_{w}=U P_{w}$ for $w \in \W_{n,k}$.  Let
$\overline{C}_w$ be the closure of the cell $C_{w}$ in Zariski
topology on on $X_{n,k}$.  Then the cohomology classes
$[\overline{C}_w]$ can be represented by variations on Schubert
polynomials and these polynomials descend to a basis of $R_{n,k}$ over
$\mathbb{Z}$ \cite[Sec. 1, Prop 3.4]{PR}.  
The Poincar\'e polynomial for $H^{*}(X_{n,k},\mathbb{Z})$ is
determined by 
\begin{equation}\label{eq:dim.fubini.mahonian}
\sum_{w \in
\W_{n,k}} q^{\codim(C_{w})} = [k]!_q \cdot \text{rev-Stir}_q(n,k),
\end{equation}
where $\text{rev-Stir}_q(n,k)$ is the polynomial obtained by reversing
the coefficients of a well-known $q$-analog of the Stirling numbers of
the second kind \cite{ER,RW,S}.

Given the impressive results due to Pawlowski and Rhoades, we call
$C_w=U P_{w}$ the \textbf{Pawlowski-Rhoades cell} or \textbf{PR
cell} indexed by $w \in \W_{n,k}$. Similarly, the \textbf{PR variety}
is denoted $\overline{C}_w$. The PR cells and PR varieties are natural
variations on the theory of Schubert cells/varieties extending to $k
\times n$ matrices, hence we believe they merit careful study of their
own.  We have used known theorems for Schubert varieties as
inspiration for conjectures and results on PR varieties.

It follows from \cite[Sec. 5]{PR} that the PR variety $\overline{C}_w$
is defined by certain bounded rank conditions.  The rank conditions
give rise to the ideal $I_{w}$ generated by the minor determinants
$\Delta_{I,J} \in \mathbb{C}[x_{11},\dotsc ,x_{kn}]$ for $I,J \in
\binom{[n]}{h}$ with $h \in [k]$ which vanish on every matrix in
$C_{w}=U P_{w}.$ The zero set of these minors is well defined on the
orbits in $\Fullrk/T$ since the right action of the diagonal matrices
just rescales each such minor.  Therefore, the spanning line
configurations in $\overline{C}_w$ can be represented by matrices in
$\Fullrk$ that vanish for every minor generating $I_{w}$.

\begin{definition}\label{def:med.roast.order}\cite[Sec. 9]{PR} The
\textbf{medium roast Fubini-Bruhat order} $(\W_{n,k},\leq)$
is defined on Fubini words by $v \leq w$ if and only if one of the
following equivalent statements is true:
\begin{enumerate}
\item $I_v \subset I_w$,
\item $\overline{C}_v \supseteq C_w$,
\item $\{(I,J) \given \Delta_{I,J}(M) =0 \ \forall M \in C_{v} \}
\subset \{(I,J) \given \Delta_{I,J}(M) =0 \ \forall M \in C_{w} \}$.
\end{enumerate}
\end{definition}

One can observe that medium roast order on Fubini words is equivalent
to Bruhat order on permutations when $n=k$.  As with Bruhat order, it
follows by definition that $v< w$ implies $\codim (C_{v})< \codim
(C_{w})$. However, some of the properties for Bruhat order on $S_{n} =
\W_{n,n}$ do not extend to all $\W_{n,k}$.  Specifically, if $v \leq
w$ in $\W_{n,k}$, then $\overline{C}_v \cap C_w \neq \emptyset$,
but the converse does not necessarily hold.  For example, using
the third condition above and the definition of pattern matrices in
\Cref{def:fubini.pattern.matrix}, one can observe that
$\overline{C}_{1323}$ contains the matrix $M_{1123} \in C_{1123}$, but
$C_{1323}$ and $C_{1123}$ are cells of the same dimension so $1323$
and $1123$ are unrelated in medium roast order.  Since $\overline{C}_v
\cap C_w \neq \emptyset$ is a weaker condition than $C_w \subseteq
\overline{C}_v$, this suggests a refinement of the medium roast
Fubini-Bruhat order, which we will denote by $\preceq$. Note that our
notation for $\preceq$ is $\leq'$ in Pawlowski and Rhoades'
notation. They use $\preceq$ for the dual order to $\leq$.

\begin{definition} For $v,w \in \W_{n,k}$, we say $C_v$
\textbf{touches} $C_w$ if $\overline{C}_v \cap C_w \neq \emptyset$,
denoted $v \rightharpoonup w$.
\end{definition}

Pawlowski and Rhoades observe in \cite[Sec. 9]{PR} that unlike the
medium roast order relations, the touching relation on Fubini words is
not transitive.  However, they showed that the transitive closure of
the touching relations is acyclic \cite[Prop. 9.2]{PR}, so the
touching relations give rise to a poset on $\W_{n,k}$ first studied but
not named in \cite{PR}.

\begin{definition} \cite[Sec. 9]{PR} The \textbf{espresso Fubini-Bruhat order}
$(\W_{n,k},\preceq)$ is defined by taking the transitive
closure of the relations of the form $v \rightharpoonup w$ if $v$ touches $w$.  
\end{definition}

Observe that for Fubini words $v, w \in \W_{n,k}$, $v \leq w$ implies
$v \preceq w$.  Thus, the medium roast order is a subposet of the
espresso order on the same set of elements.

Pawlowski and Rhoades asked for a combinatorial description of the
espresso and medium roast Fubini-Bruhat orders
\cite[Prob. 9.5]{PR}. We address this problem by giving two more sets
of defining equations for PR varieties $\overline{C}_{w}$ inside
$X_{n,k}$, see \Cref{thm:flag.minor.def} and
\Cref{thm:ranks.from.ess.set} below.  Each set is typically properly
contained in the set of all minors that vanish on the PR cell $C_{w}$,
and hence ``more efficient''.

Let $\Delta_{J}$ be the \textbf{flag minor} associated to columns in
$J$ and rows $1,2,\dotsc , |J|$.  Such minors are used historically
for the Pl\"ucker embedding of the flag variety into projective space
\cite{F}.  Note that the flag minors are invariant under the left
action of the unitriangular matrices. Hence, to determine
the vanishing/non-vanishing flag minors of $M\in C_{w}=U P_{w}$, it
suffices to consider the unique $U$-orbit representative of $M$ in
$P_{w}$.  We can partition the set of all flag minors on $k \times n$
matrices into the \textbf{sometimes, truly, and unvanishing flag
minors} for $w$, by defining the indexing sets
\begin{align*}
S_w &= \{J \in \binom{[n]}{[k]} \, | \, \exists \, A,B
\in C_w \text{ s.t. }\Delta_{J}(A) = 0, \ \Delta_{J}(B) \neq 0\},
\\ \bigskip
T_w &= \{J \in \binom{[n]}{[k]}\, | \, \Delta_{J}(M) = 0 \,\,\,
\forall \, M \in C_w\}, \text{ and}\hspace{0.55in}
\\
\bigskip
U_w &= \{J \in \binom{[n]}{[k]} \, | \,  \Delta_{J}(M)
\neq 0 \,\,\, \forall \, M \in C_w\}.
\end{align*}

\begin{theorem}\label{thm:flag.minor.def}
For any Fubini word $w \in \W_{n,k}$, the PR variety $\overline{C}_w$
is the set of spanning line configurations in $X_{n,k}$ represented by
matrices such that all  flag minors indexed by
$T_{w}$ vanish, so 
\[\overline{C}_w = \{A \in X_{n,k}\, | \, \Delta_{J}(A) = 0 \, \forall
\, J \in T_w\}.
\]
\end{theorem}

Note, the ideal $J_{w}$ generated by the flag minors $\{\Delta_{J}
\given J \in T_w\}$ is in general not the same as $I_{w}$ generated by
all vanishing minors for $C_{w}$.  For example, using the definition
and example of $P_{w}$ in \Cref{sec:background}, one can observe that
the minor $\Delta_{\{2 \},\{1 \}}=x_{21}$ is not in the ideal $J_{w}$
for $w=31123$, but it does vanish on all of $C_{w}$.  Note, both
$I_{w}$ and $J_{w}$ are radical ideals since determinants don't
factor, so they determine different affine varieties in
$\mathbb{C}^{nk}$, which agree on $X_{n,k}$.

\begin{corollary}\label{thm:Fubini.Bruhat.flag.inclusion}
For any two Fubini words $v,w\in \W_{n,k}$, we have
\begin{enumerate}
\item \label{part:medium.roast.flag} $v\leq w$ in medium roast Fubini-Bruhat order if and only if $T_v \subseteq T_{w}$, 
and
\item \label{part:espresso.flag} $v \rightharpoonup w$ if and only if $T_v \subseteq (S_w \cup T_w)$.
\end{enumerate}
\end{corollary}

Identifying vanishing flag minors of $C_w$ is more efficient
than calculating all vanishing minors of $C_w$, but still cumbersome
directly from the definition.  In fact, we can characterize the
sometimes, truly, and unvanishing flag minors via the Gale partial
order on certain multisets $\alpha_{J}(w)$ defined below.  We refer to
this as the \textbf{Alpha Test}.
These tests generalize Ehresmann's Criteria for Bruhat order in
$S_{n}$ using the Gale partial order on multisets denoted $A
\trianglelefteq B$. See \Cref{sec:background} for more details.  

\begin{definition}\label{def:alpha.vector.1} For any Fubini word $w \in
\W_{n,k}$, let $\alpha_i = \alpha_{i}(w)$ denote the position of the
initial $i$ in $w$ for each $i \in [k]$.  Call $\alpha
(w)=(\alpha_{1},\ldots, \alpha_{k})$ the \textbf{alpha vector} of $w$.
We will sometimes drop the $(w)$ when it is clear from
context. Observe that when $k=n$, the alpha vector coincides with the
notion of $w^{-1} \in S_n = \W_{n,n}$.  For $J \subset [n]$, define
the multiset 
\begin{equation}\label{def:alphaJ}
\alpha_{J}(w)=\{\alpha_{w(j)} \, | \, j \in J\}.
\end{equation}
\end{definition}

\begin{theorem}\label{thm:always.maybe.never.vanishing}(\textbf{The
Alpha Test)} Suppose $w \in \W_{n,k}$ and $J \in
\binom{[n]}{[k]}$ with $|J| = h$. Then
\begin{enumerate}
\item \label{part:sometimes} $J \in S_w$ if and only if $\{\alpha_1,\ldots,\alpha_h\} \underset{\neq}\triangleleft \alpha_J(w)$,
\item \label{part:every.time} $J \in T_w$ if and only if $\{\alpha_1,\ldots,\alpha_h\} \not\trianglelefteq \alpha_J(w)$, and
\item \label{part:unvanishing} $J \in U_w$ if and only if $\{\alpha_1,\ldots,\alpha_h\} = \alpha_J(w)$.
\end{enumerate}\end{theorem}

For example, let $w = 21231231 \in \W_{8,3}$ and $J =
\{2,6,8\}$.  Then $\alpha(w)=(\alpha_1,\alpha_2,\alpha_3) = (2, 1,4)$,
and $\alpha_{J}=\{\alpha_{w(2)},\alpha_{w(6)}, \alpha_{w(8)} \} =
\{2,1,2 \}$.  Since $\{\alpha_{1},\alpha_{2},\alpha_{3} \}=\{1,2,4
\}\not\trianglelefteq \{1,2,2\} = \alpha_J(w)$ in Gale order, we know
$J \in T_{w}$.

\medskip

\begin{corollary}\label{thm:espresso.test}
Let $v,w \in \W_{n,k}$.  Then, $v\leq w$ in medium roast Fubini-Bruhat order if and only
if for each $J \in \binom{[n]}{[k]}$ with $|J| = h \leq k$ such that
\begin{equation}\label{eq:fubini.bruhat.order}
\{\alpha_{1}(w),\ldots ,\alpha_{h}(w) \} \trianglelefteq \alpha_J(w) 
\end{equation}
we also have 
\begin{equation}\label{eq:fubini.bruhat.order.2}
\{\alpha_{1}(v),\ldots ,\alpha_{h}(v) \} \trianglelefteq
\alpha_J(v).
\end{equation}
A similar test for $v \rightharpoonup w$ holds as well based on
testing each $J$ such that $\{\alpha_{1}(w),\ldots ,\alpha_{h}(w) \} =
\alpha_J(w)$.  Therefore, if $v\leq w$ or $v\rightharpoonup w$, we
have $\{\alpha_{1}(v),\ldots ,\alpha_{h}(v) \} \trianglelefteq
\{\alpha_{1}(w),\ldots ,\alpha_{h}(w) \} $ for all $1\leq h\leq k$, 
generalizing the Ehresmann Criterion.

\end{corollary}

In \Cref{sec:background}, we briefly review our notation and key
concepts from the literature.  In \Cref{sec:proofs}, we indicate some
of the lemmas needed to prove \Cref{thm:flag.minor.def} and its
corollaries.  In \Cref{sec:covering}, we identify certain families of
covering relations and use them to define the decaf Fubini-Bruhat
order. We also state an analog of the Lifting Property of Bruhat
order.  In \Cref{sec:essential}, we generalize Fulton's essential set
for permutations to Fubini words and show this set gives the unique
minimal set of rank conditions defining a PR variety, see
\Cref{thm:ess.set.with.rank}.

\section{Background}\label{sec:background}

For a positive integer $n$, let $[n]$ denote the set $\{1,2,\dots , n
\}$.  Generalizing the notation for binomial coefficients, we let
$\binom{[n] }{k}$ denote all size $k$ subsets of $[n]$ and $\binom{[n]
}{[k]} =\bigcup_{h=1}^{k} \binom{[n] }{h}$.  The \textbf{Gale order}
on $\binom{[n] }{k}$ is given by $\{a_{1}< \dots < a_{k}\}
\trianglelefteq \{b_{1}< \dots < b_{k}\}$ if and only if $a_{i}\leq
b_{i}$ for all $i \in [k]$ \cite{Gale}.  Gale order can easily be
extended to multisets of positive integers of the same size.

Let $S_{n}$ denote the symmetric group on $[n]$ thought of as
bijections $w:[n]\to [n]$.  As usual, write a permutation $w$ in
\textbf{one-line notation} as $w=w_{1}\cdots w_{n}$.  Let $t_{ij}$ be
the transposition interchanging $i$ and $j$, and let $s_{i}$ denote
the simple transposition interchanging $i$ and $i+1$.  The permutation
$t_{ij}w$ is obtained from the one-line notation for $w$ by
interchanging the values $i$ and $j$, while right multiplication
$wt_{ij}$ interchanges the values $w_{i}$ and $w_{j}$.  The
permutation matrix $M_{w}$ for $w \in S_{n}$ is the $n \times n$
matrix with a 1 in position $(w_{j},j)$ for all $j \in [n]$ and 0's
elsewhere.  Permutation multiplication agrees with matrix
multiplication: $u=vw$ if and only if $M_{u}=M_{v}M_{w}$.  Permutation
multiplication extends to Fubini words if the corresponding matrices
have the correct size.

Schubert varieties $X_{w}$ for $w \in S_{n}$ in the flag variety
$GL_{n}/B_{+}^{(n)}$ are defined via bounded rank conditions on matrices
coming from the associated permutation matrices \cite{F}.  The
\textbf{Bruhat order} on $S_{n}$ is defined by reverse inclusion on
Schubert varieties: $v\leq w \iff X_{w}\subset X_{v}$.  This poset can
be characterized as the transitive closure of the relation $w\leq
t_{ij}w$ provided $i<j$ and $i$ appears to the left of $j$ in the
online notation for $w$ \cite{BB}.  The covering relations are given
by the set of edges $w\leq t_{ij}w$ such that $t_{ij}w$ has exactly
one more inversion than $w$.  Ehresmann characterized Bruhat order on
$S_{n}$ in terms of Gale order, decades prior to Gale or Bruhat's
work, by the \textbf{Ehresmann Criterion} \cite{E}
\begin{equation}\label{eq:bruhat}
v\leq w \iff \{v_{1},v_{2},\dots , v_{i} \} \trianglelefteq
\{w_{1},w_{2},\dots , w_{i} \} \ \forall i \in [n].  
\end{equation}
Suppose $v \leq w$ in Bruhat order on $S_{n}$, $i \in [n-1]$ and $i+1$
precedes $i$ in both $v$ and $w$. Then, the \textbf{Lifting Property of Bruhat
order} \cite[Prop. 2.2.7]{BB} implies that $s_i v \leq s_i w$.

\begin{definition}\label{def:diagram.def}
The \textbf{Rothe diagram} of a permutation $w \in S_{n}$ is the
subset of $[n]\times [n]$ in matrix coordinates given by $D(w) =
\{(w_j,i) \, | \, i<j \text{ and } w_i > w_j\}.$ Define the
\textbf{essential set} of $w$, denoted $Ess(w)$, to include all $(i,j)
\in D(w)$ such that $(i+1,j),(i,j+1) \not \in D(w)$.
\end{definition}

The Rothe diagrams are used extensively in the theory of Schubert
varieties.  In particular, Fulton showed that the rank conditions
coming from the coordinates $(i,j) \in Ess(w)$ determine the unique
minimal set of bounded rank equations defining the Schubert variety $X_{w}$
\cite{FES}. Eriksson-Linusson showed that the average size of the
essential set is $n^{2}/36$ for $w \in S_{n}$ \cite{EL}.

Much of the notation for permutations defined above has an analog for
Fubini words.  For $w = w_{1}\cdots w_{n} \in \W_{n,k}$, let $M_{w}$ be
the matrix obtained from the $k \times n$ all zeros matrix by setting
the $(w_{j},j)$ entry to be 1 for all $j \in [n]$. 
Note that $M_{w}$ has exactly one 1 in each column and at least one 1 in
each row, but it may have many 1's in any row.  Recall from
\Cref{def:alpha.vector.1} that $\alpha_{i}(w)=\alpha_{i}$ is the
position of the first letter $i$ in $w$ for $i \in [k]$.

\begin{definition}\label{def:initial.positions} \cite[\S 3]{PR} For a
word $w \in \W_{n,k}$, the \textbf{initial positions} of $w$ are the
set $\text{in}(w) = \{\alpha_1,\ldots,\alpha_k\}$. A \textbf{redundant
position} of $w$ is any position that is not initial. An
\textbf{initial letter} is a letter appearing in an initial position,
and a \textbf{redundant letter} is a letter appearing in a redundant
position.
\end{definition}

\begin{definition}\label{def:initial.permutation} \cite[\S 3]{PR} For
$w \in \W_{n,k}$, the \textbf{initial permutation}, $\pi(w) \in S_{k}$,
is obtained from $w$ by deleting the redundant letters from the
one-line notation. 
\end{definition}

\begin{definition}\label{def:fubini.pattern.matrix} \cite[\S 3]{PR}
For $w = w_{1}\cdots w_{n} \in \W_{n,k}$, the \textbf{pattern matrix}
$P_w$ is a $k \times n$ matrix with entries $0$, $1$, or
$\star$. Obtain $P_w$ by starting with $M_w$ and replacing the $0$ by
a $\star$ in each position $(w_i,j)$ such that $i \in \text{in}(w)$,
$i < \alpha_{w(j)}$, $1\leq j\leq n$, and either $j \in in(w)$ and
$w_i<w_j$, or $j \not\in in(w)$.

A matrix is said to \textbf{fit the pattern of $w$} if that matrix can
be obtained by replacing the $\star$'s in the pattern matrix of $w$
with complex numbers.  We will abuse notation and consider $P_{w}$
both as a $k \times n$ matrix with entries in $\{0,1,\star \}$ and as
the set of all matrices fitting the pattern of $w$. 
\end{definition}

\begin{definition}\cite[Eq. (3.6)]{PR} The \textbf{dimension} of $w
\in \W_{n,k}$, denoted $\text{dim}(w)$, is the number $\star$'s in its
pattern matrix $P_w$.
\end{definition}

\begin{example}\label{example:fubini.pattern.matrix}
The pattern matrices of $v=31422$ and $w=31424$ in $\W_{5,4}$ are
\[
P_{31422} =
\begin{pmatrix}
0 & 1 & \star & \star &  \star  \\
0 & 0 &     0 &     1 &  1  \\
1 &0  &     \star &     0 & \star \\
0 & 0 &     1 &     0 &  \star
\end{pmatrix}
\text{  and  }
P_{31424} =
\begin{pmatrix}
0 & 1 & \star & \star &  \star  \\
0 & 0 &     0 &     1 &  0  \\
1 &0  &     \star &     0 & \star \\
0 & 0 &     1 &     0 &  1
\end{pmatrix}.
\]
Therefore, $\text{dim}(31422)=6$ and $\text{dim}(31424)=5$.  
\end{example}

If $w \in \W_{n,k}$, then the dimension of the PR cell $C_{w}$ is
$\mathrm{dim}(w)+\binom{k }{2}$.  The unique largest dimensional cell
in $X_{n,k}$ is $C_{123\cdots kk^{n-k}}$ and $\text{dim}(12\cdots k
k^{n-k})=\binom{k}{2}+(n-1)(k-1)$.  Hence,
$X_{n,k}=\overline{C}_{12\cdots k k^{n-k}}$ has dimension
$n(k-1)=2\binom{k}{2}+(n-1)(k-1)$ and $12\cdots k k^{n-k}$ is the
unique minimal element in all three Fubini-Bruhat orders.  Since
Fubini words are in bijection with ordered set partitions, the
dimension generating function gives a natural $q$-analog of the
Stirling numbers of the second kind $ \sum_{w \in \W_{n,k}}
q^{\mathrm{dim}(w)} = [k]!_q \cdot \text{Stir}_q(n,k).$ Reversing the
coefficients in this generating function gives
\eqref{eq:dim.fubini.mahonian}.


\section{Outlines of Proofs}\label{sec:proofs}

We outline the proofs of \Cref{thm:flag.minor.def} and
\Cref{thm:always.maybe.never.vanishing}. These statements form the
basis from which the covering relations and other Fubini-Bruhat order
properties can be proved.

\begin{lemma}\label{lem:word.from.minors}
Given $A \in \Fullrk$, the \textit{projective coordinates} $
P(A)=(\Delta_{J}(A)\given J \in \binom{[n] }{[k]})$ determine both the
unique $w \in \W_{n,k}$ such that $A \in UP_{w}T^{(n)}$ and $A' \in
P_{w}$ such that $A\in UA'$.
\end{lemma}

\begin{corollary}\label{cor:Fubini.Bruhat.every.time.vanishing}
The set $T_{w}$ of truly vanishing flag minors on the PR cell $C_w$
determines $w \in \W_{n,k}$, and therefore the rank conditions
defining $\overline{C}_{w}$ as a subset of $X_{n,k}$. 
\end{corollary}

\Cref{cor:Fubini.Bruhat.every.time.vanishing} says there is
enough information in the set $T_{w}$ to recover $w$.  To make the
relationship between $T_{w}$ and $\overline{C}_{w}$ precise, we
observe several relations among minors that hold specifically on PR
cells and spanning line configurations.

\begin{lemma}\label{lem:sub.col.dependence}
Suppose $w \in \W_{n,k}$ is a Fubini word, $J \subset [n]$, and $1\leq
h\leq k$.  Let $\text{rank}_w^{(h)}(J)$ be the largest value $r$ such
that there exist subsets $I \subset [h]$ and $J' \subset J$ such that
$r=|I|=|J'|$ and $\Delta_{I,J'}(A) \neq 0$ for some $A \in C_{w}$.
The following conditions are equivalent.

\begin{enumerate}
\item We have $\text{rank}_w^{(h)}(J) < |J|$.

\item For every $I \subseteq [h]$ such that $|I|=|J|$, the
$(I,J)$-minor vanishes on $C_{w}$.

\item For all subsets $K \in \binom{[n] }{h}$ such that $J\subset
K$, we have  $K \in T_w$.
\end{enumerate}
\end{lemma}

\begin{corollary}\label{cor:big.minor}
Suppose $w \in \W_{n,k}$ is a Fubini word, $I\subseteq [k]$
and $J\subseteq [n]$ are sets of the same size, and $h=\max(I)$. If
the $(I,J)$-minor vanishes on $C_{w}$, then at least one of
the following hold.

\begin{enumerate}
\item \label{part:upper.minors} For every $j \in J$, the $(I \setminus \{h\},J \setminus \{j \})$-minor vanishes on $C_{w}$.
\item \label{part:sub.cols} For all subsets $K$ such that $J\subseteq
 K\in \binom{[n] }{h}$,  we have $K \in T_w$.
\end{enumerate}
\end{corollary}

 \Cref{cor:big.minor} follows from \Cref{lem:sub.col.dependence}.
\Cref{thm:flag.minor.def} follows by induction on the number of rows
of a minor of $C_w$ using \Cref{cor:big.minor}, and by
\Cref{lem:sub.col.dependence}.

\begin{lemma}\label{lem:one.matrices}
Suppose $w \in \W_{n,k}$ is a Fubini word and $J \in \binom{[n]}{[k]}$
with $h = |J|$. Then, $J \in U_w$ if and only if the submatrix
$M_{w}[[h],J]$ is a permutation matrix.
\end{lemma} 

\begin{lemma}\label{lem:make.flag}
Let $w \in \W_{n,k}$, $I\subseteq [k]$ and $J\in \binom{[n]}{[k]}$
such that $|I|=|J|$ and $\Delta_{I,J}(A)=0$ for all $A$ in the PR cell
$C_{w}$. Then $(H,J)$ indexes a vanishing minor on $C_{w}$ for
any $H$ such that $|H|=|I|$ and $H \leq_L I$ in lex order. In
particular, $\Delta_{[|I|],J} $ is a vanishing flag minor on $C_{w}$,
so $J \in T_w$.
\end{lemma}

Lemmas~\ref{lem:one.matrices} and~\ref{lem:make.flag}, together with
the earlier lemmas can be used to prove Corollary
\ref{thm:Fubini.Bruhat.flag.inclusion}. Corollary
\ref{thm:Fubini.Bruhat.flag.inclusion} and Lemma
\ref{lem:one.matrices} imply Theorem
\ref{thm:always.maybe.never.vanishing}.

\section{Covering Relations and the Decaf Order}\label{sec:covering}

The following rules describe some families of covering relations for
the medium roast and espresso Fubini-Bruhat orders, giving a partial
answer to Problem 9.5 in \cite{PR}. The Transposition Rule and the
Pushback Rule allow us to define the decaf Fubini-Bruhat order, the
only ranked Fubini-Bruhat order. We also discuss a generalization of
the Lifting Property from Bruhat order.

We start with two observations on covering relations that follow from
the definition of medium roast order, pattern matrices, and
\Cref{thm:Fubini.Bruhat.flag.inclusion}.  Let $w = w_1 \cdots w_n \in
\W_{n,k}$ with initial permutation $\pi (w)=\pi_{1}\cdots \pi_{k}$.
\begin{enumerate}
\item \textbf{The Transposition Rule.}  For $1\leq i<j\leq k$, we have
$w < t_{ij} w$ in medium roast Fubini-Bruhat order if and only if
$\alpha_i(w) < \alpha_j(w)$. In particular, $t_{ij} w$ covers $w$ in
medium roast Fubini-Bruhat order if and only if $\pi(t_{ij} w)$ covers
$\pi(w)$ in Bruhat order on $S_k$.

\item
\textbf{The Pushback Rule.} Suppose $w_j = \pi_i$ is a redundant
letter in $w$ for $i \in [k-1]$ and $j \in [n]$.  Let $v$ be the
Fubini word obtained from $w$ by replacing $w_{j}$ by $\pi_{i+1}$.
Then, $w$ covers $v$ in medium roast Fubini-Bruhat order.  See
\Cref{example:fubini.pattern.matrix} for an example of $v<w$
satisfying the pushback covering relation.
\end{enumerate}

\begin{definition} The \textbf{decaf Fubini-Bruhat order} on
$\W_{n,k}$ is the transitive closure of the covering relations given
by the Transposition Rule and the Pushback Rule.
\end{definition}

The decaf order has many nice properties. It is the product of Bruhat
order for $S_{k}$ and the poset determined by pushbacks on the subset
$\{w\in \W_{n,k}\given \pi (w)=id \}$.  The decaf order is a ranked
poset on $\W_{n,k}$, and its rank generating function is the same as
the Poincar\'e polynomial in \eqref{eq:dim.fubini.mahonian}.  The
medium roast and espresso orders are not ranked posets in general. For
$n\geq 5$ and most values of $k$, there are covering relations in the
medium roast Fubini-Bruhat order $(\W_{n,k},\leq)$ with a dimension
difference of 2 or more, causing the medium roast Fubini-Bruhat order
to be unranked in general. For example, in $\W_{5,4}$, 44312 covers
41321, but 44312 has dimension 1, and 41321 has dimension 3.

\begin{theorem}\label{thm:the.superpushback.rule}
\textbf{The Superpushback Rule.} Suppose $w \in \W_{n,k}$, $i \in
[k-1]$, and $j \in [n]$ such that $w_j = \pi_i$ is a redundant letter
in $w$. If $i+p \leq k$ and $v$ is obtained from $w$ by replacing
$w_{j}$ by $\pi_{i+p}(w)$, then $v \rightharpoonup
w$ and this is a covering relation in both espresso and medium roast orders.\end{theorem}

\begin{theorem}\label{thm:the.lifting.property}
\textbf{The Lifting Property.} Suppose $v,w \in \W_{n,k}$, $i \in
[k-1]$, $\alpha_{i+1}(v) < \alpha_i(v)$, and $\alpha_{i+1}(w) <
\alpha_i(w)$.  If $v \leq w$ in medium roast Fubini-Bruhat order, then
$s_i v \leq s_i w$.  Furthermore, if $v \rightharpoonup w$, then $s_i
v \rightharpoonup s_i w$.
\end{theorem}

\section{Essential Sets}\label{sec:essential}

We extend the notion of a Rothe diagram from \Cref{def:diagram.def} to
Fubini words. This allows us to define the essential set for a Fubini
word.  We then show the essential set determines a minimal set of rank
equations on the corresponding PR variety, generalizing Fulton's
essential set for permutations and Schubert varieties \cite{FES}.
This leads to an essential set characterization of $v\leq w$ in medium
roast order.

\begin{definition}\label{def:conv} \cite{PR}
A Fubini word $w \in \W_{n,k}$ is called \textbf{convex} if $h<j$ and
$w_h=w_j$ implies that $w_i=w_j$ for every $i$ such that $h<i<j$. Then
the \textbf{convexification} of $w$, denoted by $\text{conv}(w)$, is
the unique convex word such that $\pi(\text{conv}(w)) = \pi(w)$ and
the content of $w$ and $\text{conv}(w)$ are the same as multisets.
The \textbf{standardization} of $w$, denoted $\text{std}(w) \in S_n$,
is obtained by replacing the $n-k$ redundant letters of $w$ with
$k+1,k+2,\ldots,n$ from left to right.
\end{definition}

 Deduce from \Cref{def:conv} that two Fubini words $v,w \in \W_{n,k}$
have the same convexification, $\text{conv}(v)=\text{conv}(w)$, if and
only if $\pi (v)=\pi (w)$ and they have the same multiset of
letters.

\begin{definition}\label{def:diagram}
Given Fubini word $w \in \W_{n,k}$, define the \textbf{diagram} of $w$ to be $D(\text{std}(\text{conv}(w)))$.
\end{definition}

\begin{figure}[h]
  \centering
  \includegraphics[width=0.2\textwidth]{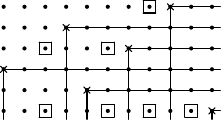}
  \caption{Diagram of $44253136541$ with cells in the essential set boxed.}
  \label{fig:44253136541.diagram}
\end{figure}

One can observe that $D(\text{std}(\text{conv}(w))) \subset [k] \times
[n]$, as none of the bottom $n-k$ rows will contribute any elements to
$D(\text{std}(\text{conv}(w)))$.  Thus, the diagram of a Fubini word
in $\W_{n,k}$ can be drawn as a $k \times n$ grid of dots.  For
example, the convexification of $w=44253136541 \in \W_{11,6}$ is
$44425533116$, and
$\text{std}(44425533116)=[4,7,8,2,5,9,3,10,1,11,6]$.  So the diagram
for $w$ is $D([4,7,8,2,5,9,3,10,1,11,6])$.   See Figure~\ref{fig:44253136541.diagram}.

In analogy with the alpha vector, define the \textbf{beta vector}
$\beta (w)=(\beta_{1}(w),\dotsc , \beta_{k}(w))$ for $w \in W_{n,k}$
by $\beta_i(w) =\beta_{i}= \{j \in [n] \given w_{j} \in
\{\pi_{1},\dotsc , \pi_{i} \} \}$ where $\pi (w)=(\pi_{1}, \dotsc ,
\pi_{k}) \in S_{k}$ is the initial permutation.  Note that $\beta_{1}
\subset \cdots \subset \beta_{k}$.  For example, if $w = 12123123 \in
\W_{8,3}$, we observe $\beta_1 = \{1,3,6\}$, $\beta_2 = \{1,2,3,4,6,7
\}$, and $\beta_3\, =[8]$.

Given any Fubini word $w \in \W_{n,k}$, define its \textbf{rank
function} to be the map $r_w : [k] \times [k] \rightarrow
\mathbb{Z}_{\geq 0}$ that sends $(h,i)$ to the maximum value of the
rank of the submatrix $A[[h],\beta_{i}]$ over all $A \in C_{w}$.
This function can be determined directly from the Fubini word $w$ as
with permutations, but the statement is more complicated so we have
omitted it for brevity. From the pattern matrix definition, one can
observe that the jumps in the rank functions of matrices in a PR
variety are determined by the sets in the beta vector.

\begin{definition}\label{defn:essential.set}
Given any Fubini word $w \in \W_{n,k}$, define the \textbf{ranked
essential set} of $w$ to be
\[
Ess^{*}(w)=\{(h, \beta_i, r) \, | \, (h,|\beta_{i}|) \in
Ess(\text{std}(\text{conv}(w))), r=r_w(h,i) \}.
\]
\end{definition}

\begin{theorem}\label{thm:ranks.from.ess.set} A matrix $A \in \Fullrk$
is in the PR variety $\overline{C}_w$ if and only if the rank of the
top $h$ rows of $A$ in the columns $\beta_i(w)$ is at most $r$ for
each $(h,\beta_i(w),r) \in Ess^{*}(w)$, and no smaller set of rank
conditions will suffice.
\end{theorem}

\begin{corollary}\label{thm:ess.set.with.rank} Let $v, w \in \W_{n,k}$. Then $v \leq w$ if and only if for every $(m,\beta_j(v),s) \in Ess^{*}(v)$, there exists an $(h,\beta_i(w),r) \in Ess^{*}(w)$ such that
$\text{max}(0,m-h) + | \beta_j(v)
\setminus \beta_i(w) |  \leq s-r.$
\end{corollary}

Bj\"orner-Brenti gave an improvement on the Ehresmann Criterion for
Bruhat order on permutations in \cite{BB2}.  Similar improvements on
the Alpha Test for medium and espresso orders exist as well.  Such
improvements also lead to a reduction in the number of equations
necessary to define a PR variety. In recent work, Gao-Yong found a
minimal number of equations defining a Schubert variety in the flag
variety \cite{Gao-Yong}.  Thus, it would be interesting to consider
the following problem.

\begin{open}
Identify a minimal set of equations defining a PR variety.
\end{open}

\acknowledgements{ This paper grew out of the Ph.D. thesis of the
second author \cite{Billey-Ryan,Ryan.2022}.  Further details and more background
can be found there.  Many thanks to Brendan Pawlowski, Brendan
Rhoades, Jordan Weaver, and the referees for insightful suggestions on
this project.}

\printbibliography

\end{document}